\documentclass[12pt]{article}
\usepackage[utf8]{inputenc}
\usepackage{amsmath}
\usepackage{amsfonts}
\usepackage{amssymb}
\usepackage{amsthm}
\usepackage{tipa}
\usepackage{relsize}
\usepackage{array}
\usepackage{multirow}
 \usepackage{csquotes} 
\usepackage{tabu}
\usepackage{soul}
\usepackage{graphicx}
\usepackage{epigraph}
\usepackage{float}
\usepackage[english]{babel}
\usepackage[usenames, dvipsnames]{color}
\usepackage{fancyhdr}
\usepackage[Sonny]{fncychap}
\usepackage{tikz-cd}
\usepackage{tkz-euclide}
\usepackage{sseq}
\usepackage{newfloat}
\usepackage{amsmath,mathtools}
\usepackage{pgfplots}
\usepackage{geometry}
 \usepackage{sectsty}
\sectionfont{\fontsize{14}{15}\selectfont}

\usepackage{endnotes}

\usepackage[all]{xy}
\DeclareFloatingEnvironment[fileext=dia]{diagram}
\theoremstyle{definition}

\usepackage[most]{tcolorbox}

\usepackage{tabularx}
\usepackage{imakeidx}
\usepackage{hyperref}
\hypersetup{colorlinks={true},linkcolor={blue},citecolor={blue},urlcolor={blue}, linktoc=all}

\begin{document}
\renewcommand{\refname}{Notes}

\title{\textbf{Sexagesimal calculations in ancient Sumer}}

\author{Kazuo MUROI}

\maketitle

\begin{abstract}
 This article  discusses the reasons for the choice of the sexagesimal system by ancient Sumerians. It is shown that Sumerians chose this specific numeral system based on logical and practical reasons which enabled them to deal with big numbers easily and even perform the  multiplications and divisions in this system. I shall also discuss how the Sumerians calculated the area of a large field and measured a large quantity of barley according to their seemingly complicated but really systematic methods. %Moreover, I would like to show that it is not recommended transcribing the sexagesimal numbers occurred in their documents into the corresponding decimal numbers, which   most Sumerologists usually do, because such a transcription obscures the calculation process performed by the Sumerians and may mislead us about the sexagesimal structure of their mathematics.
\end{abstract}

\section{The advent of the sexagesimal notation}
By the 26th century BCE when historic times began, the Sumerians had invented the sexagesimal system based on their decimal system. At the latest in the 20th century BCE, they had completed the  sexagesimal place value notation, which may have been the forerunner of our decimal system now in use, though having lacked symbols both for a “sexagesimal point” and for the number zero. In the measurement of angles in degrees, minutes, and seconds and of time in hours, minutes, and seconds we still use the sexagesimal notation, which is a legacy of Sumer. This fact may show us not only that old habits die hard but also that the sexagesimal system has certain advantages over the decimal system which only resulted from the number of fingers of a person or was a result of chance in a sense.

In Sumerian society it was necessary to count and record a large amount of farm and marine products and then to distribute the products to workers or soldiers according to some rule. This process necessitated the Sumerians developing a useful and efficient numeral system, that is, the sexagesimal numeral system. Thus, their system has two noticeable characteristics. First, it fits the expressing large numbers which occur in administrative documents. Generally speaking, a number of three figures composed of $ 60^{2} $, 60 and 1 seems to have been sufficient for their task. Secondly, in this system division is easy to carry out as compared with the decimal system because the base 60 is divisible by 2, 3, 4, 5, 6, 10, 12, 15, 20, and 30. Therefore, we can say that the sexagesimal numeral system was not a product of chance, but a deliberately devised system suited for economic activities in Sumer.

Through the two tablets written in the 26th century BCE, we can assume that the Sumerian scribes were able to perform multiplication and division easily.

\begin{figure}[h]
	\centering
	\includegraphics[scale=1]{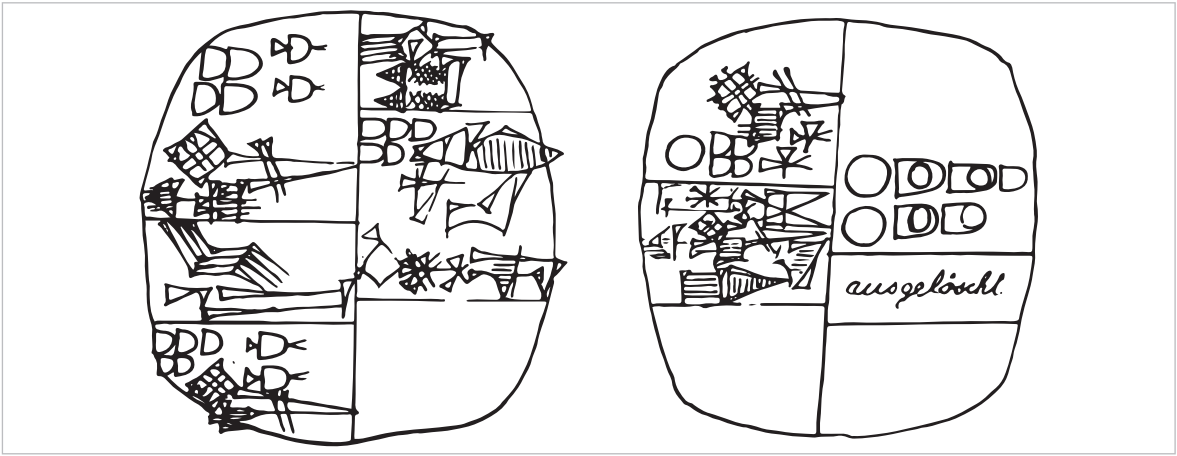}
	\caption{The obverse and reverse of WF No.2}
	\label{Figure1}
\end{figure}

\subsubsection*{WF No.2 (administrative tablet)$^{\text{\cite{AD}}}$} %\endnote{A. Deimel, Wirtschaftstexte aus Fara (= WF), 1924, p.1. Deimel erroneously interprets the number 9660 as the number of donkeys. See p.10. The translation is mine.}}
4 yellow-donkeys for bundles of flax of {\fontfamily{qpl}\selectfont Lumma} (personal name). 5 yellow-donkeys for bundles of flax of {\fontfamily{qpl}\selectfont Billum}(?). 5 (yellow-donkeys for bundles of flax) of {\fontfamily{qpl}\selectfont Lunumudae}. 2,41,0 ($ =2 \times 60^2+41\times60 = 9660= 14\times690 $). (In total) 14 (yellow)-donkeys for bundles of flax. As a consequence of the transport of bundles of flax by breeding female yellow-donkeys (see Figure \ref{Figure1}). The scribe of this tablet must have calculated $ 14\times(11,30) = 2,41,0 $, which is the total number of bundles transported by 14 donkeys.

\begin{figure}[h]
	\centering
	\includegraphics[scale=1]{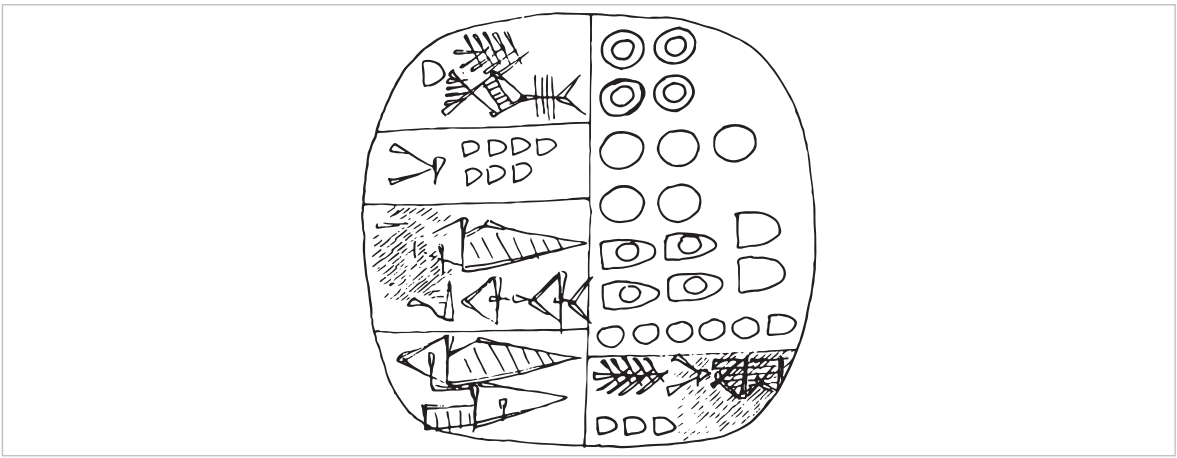}
	\caption{The obverse of TSŠ No.50}
	\label{Figure2a}
\end{figure}

\subsubsection*{TSŠ No.50 (school tablet)$^{\text{\cite{RJ}}}$%\endnote{R. Jestin, Tablettes sumériennes de Šuruppak, 1937, no.50. For the transliteration and translation, see my paper:\\
	%K. Muroi, The Origin of the Mystical Number Seven in Mesopotamian Culture: Division by Seven in the sexagesimal Number System, 2014, \hyperlink{https://arxiv.org/abs/1407.6246}{arXiv: 1407.6247 [math.HO]}.}
}
\begin{figure}[h]
	\centering
	\includegraphics[scale=1]{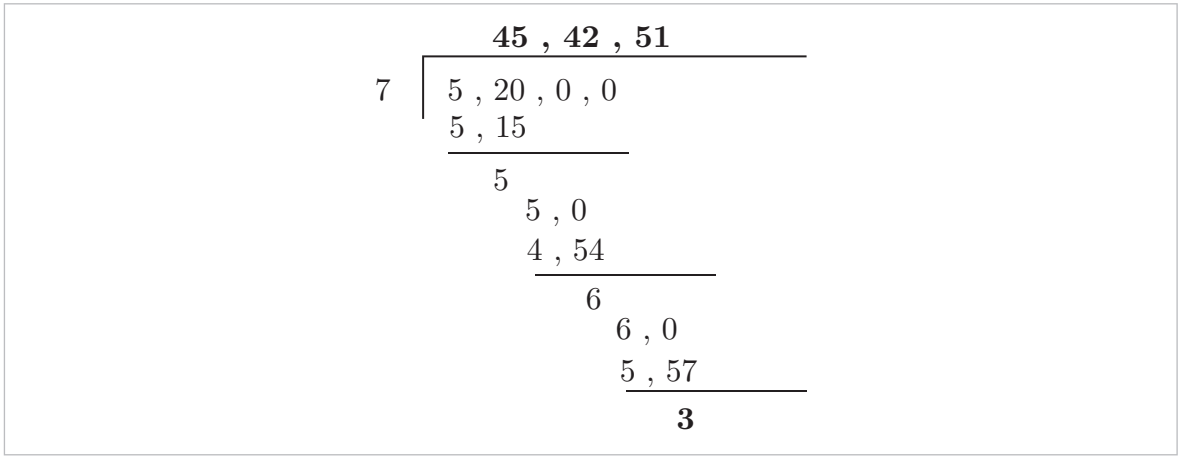}
	\caption{The assumed division process}
	\label{Figure2b}
\end{figure}

1 {\fontfamily{qpl}\selectfont gur$_{7}$} (=5,20,0,0 {\fontfamily{qpl}\selectfont sìla}, where 1 {\fontfamily{qpl}\selectfont sìla} $\approx$ 1 liter) of barley. 7 {\fontfamily{qpl}\selectfont sìla} (of    barley each) one man received. (The number of) the men are
45,42,51 ($  = 45\times 60^{2}+42\times 60+51 = 164571 $).
3 {\fontfamily{qpl}\selectfont sìla} of barley is repaid.
It is evident that the Sumerian scribe of this tablet has carried out the following division:
$ 5,20,0,0 = 7\times (45,42,51)+3 $ (see Figures \ref{Figure2a} and \ref{Figure2b}).

\section{Numerals}
The Sumerian word {\fontfamily{qpl}\selectfont li-mu-um} “one thousand” was borrowed from the Akkadian {\fontfamily{qpl}\selectfont \textit{līmu}}$^{\text{\cite{CAD9}}}$%\endnote{The Assyrian Dictionary of the Oriental Institute of the University of Chicago (= CAD), vol.9, L, 1973, p.197.}
, and also the Sumerian word me-at “one hundred” is a loanword from the Akkadian {\fontfamily{qpl}\selectfont \textit{me’atum}}$^{\text{\cite{CAD10}}}$.%\endnote{CAD, vol.10, M part two, 1977, p.1.}
Moreover, there was no single word in Sumer that designates a myriad “ten thousand” as well as in Babylonia proper. From these facts we may conclude that the Sumerians did not develop any elaborate decimal numeral system for large numbers and invented the sexagesimal system instead. The basic numerals of the Sumerian sexagesimal system, which must have been sufficient for their accounting, are as follows$^{\text{\cite{MLT}}}$%\endnote{M. L. Thomsen, The Sumerian Language, 1984, §139.}
: \\
{\fontfamily{qpl}\selectfont aš} or {\fontfamily{qpl}\selectfont diš} = 1\\
{\fontfamily{qpl}\selectfont min} = 2\\
{\fontfamily{qpl}\selectfont eš$_{5}$} = 3\\
{\fontfamily{qpl}\selectfont limmu} = 4\\
{\fontfamily{qpl}\selectfont iá} = 5\\
{\fontfamily{qpl}\selectfont àš} = 6 ($ < $ {\fontfamily{qpl}\selectfont iá+aš})\\
{\fontfamily{qpl}\selectfont imin} = 7 ($<$ {\fontfamily{qpl}\selectfont iá+min})\\
{\fontfamily{qpl}\selectfont ussu} = 8\\
{\fontfamily{qpl}\selectfont ilimmu} = 9 ($<$ {\fontfamily{qpl}\selectfont iá+limmu})\\
{\fontfamily{qpl}\selectfont u} =10\\
{\fontfamily{qpl}\selectfont niš} = 20\\
{\fontfamily{qpl}\selectfont ùšu} = 30\\
{\fontfamily{qpl}\selectfont nimin} = 40 ($<$ {\fontfamily{qpl}\selectfont niš+min})\\
{\fontfamily{qpl}\selectfont ninnu} = 50 ($<$ {\fontfamily{qpl}\selectfont nimin+u})\\
{\fontfamily{qpl}\selectfont géš} = 1,0 (= 60)\\
{\fontfamily{qpl}\selectfont šár} = 1,0,0 (= $60^{2}$)\\
{\fontfamily{qpl}\selectfont šár-gal} “big {\fontfamily{qpl}\selectfont šár}” or {\fontfamily{qpl}\selectfont šár}$\times${\fontfamily{qpl}\selectfont géš} “sixty {\fontfamily{qpl}\selectfont šár}” = 1,0,0,0 (= $ 60^{3} $) 

In combination with these numerals, the Sumerian scribes seem to have been able to pronounce any whole number smaller than $ 60^{4} $, for example:\\
{\fontfamily{qpl}\selectfont gešta-u} ($<$ {\fontfamily{qpl}\selectfont géš-ta u}) “ten after sixty” = 1,10 (= 70),\\
{\fontfamily{qpl}\selectfont géš-u} “ten sixties” = 10,0 (= 600),\\
{\fontfamily{qpl}\selectfont šár-min} (2,0,0) {\fontfamily{qpl}\selectfont géš-u-limmu} (40,0) {\fontfamily{qpl}\selectfont géš} (1,0) = 2,41,0 (= 9660) (see Figure \ref{Figure1}).

Although the largest numeral attested so far is $ 60^{4} $, it is obvious that this number word was not used in everyday life:\\
{\fontfamily{qpl}\selectfont šár-gal šu-nu-tag-ga} “big {\fontfamily{qpl}\selectfont šár} that is not touched”.$^{\text{\cite{HVH}}}$%\endnote{H. V. Hilprecht, Mathematical,  Metrological and Chronological Tablets from the Temple Library of Nippur, 1906, text 29, obverse col.4, line 9.}
\\
As to {\fontfamily{qpl}\selectfont šár} and {\fontfamily{qpl}\selectfont šár}$\times${\fontfamily{qpl}\selectfont géš}, they were borrowed into Akkadian as {\fontfamily{qpl}\selectfont \textit{šār}} ($ 60^{2} $) and {\fontfamily{qpl}\selectfont \textit{šuššār}} ($ 60^{3} $) respectively. 

A few words about the origin of the Akkadian numeral {\fontfamily{qpl}\selectfont \textit{nēr}} ``10,0'', which was adopted by the Greeks as neros, should be given. It seems that the Babylonians did not use the Sumerian numeral {\fontfamily{qpl}\selectfont géš-u} ``ten-sixty'' as it was. Instead they coined the {\fontfamily{qpl}\selectfont \textit{nēr}} because they had already used the Akkadian word {\fontfamily{qpl}\selectfont \textit{giššu}} ``hip''. Since {\fontfamily{qpl}\selectfont giš-(giš) = \textit{nīrum}} ``yoke'', they used the new word  {\fontfamily{qpl}\selectfont \textit{nēr}} for ``10,0'' and  ignored its original meaning.

\section{Number signs}
\begin{figure}[h]
	\centering
	\includegraphics[scale=1]{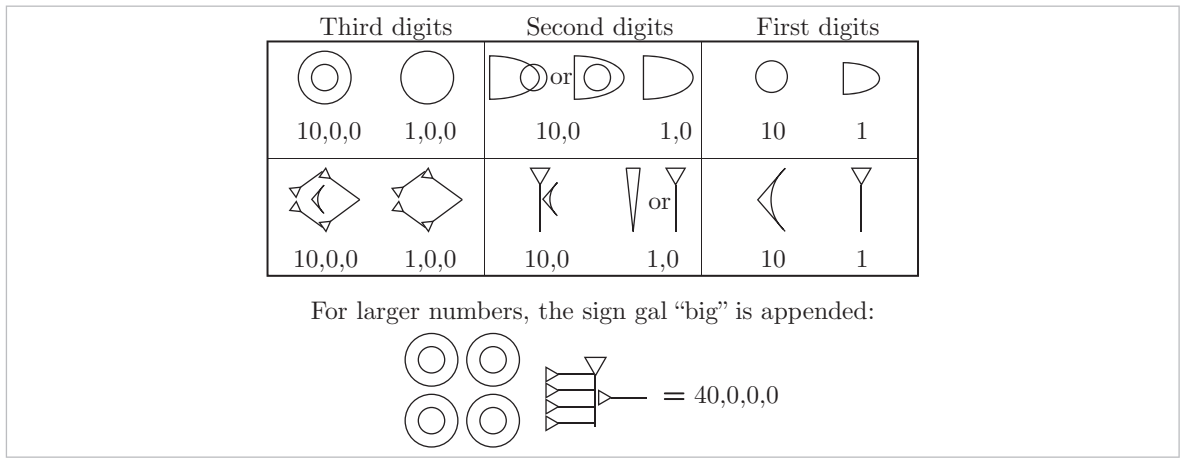}
	\caption{Number signs}
	\label{Figure3}
\end{figure}

Although the number signs used in sexagesimal notation vary in shape depending on when and where they had been written down, we would recognize two basic types of them, that is, round-shaped number signs and wedge-shaped ones (see Figure \ref{Figure3}). The former was used mainly in Early Dynastic Period III (EDIII 2600-2340 BCE) and the latter was used from the 24th century BCE onward. Naturally, there are many cuneiform tablets in which both signs occur side by side, and also some tablets of EDIII in which only wedge-shaped number signs occur.

With these number signs the Sumerian scribes recorded their results of sexagesimal calculations just as the Japanese who still sometimes write down the results of calculations, especially the sum of money, with Chinese numerals instead of Hindu-Arabic numerals.

\section{Regular numbers}
It is most probable that the Sumerians already knew from experience that the numbers 2,3, and 5 are regular to the base 60, because they especially called the numbers {\fontfamily{qpl}\selectfont a-rá-gub-ba} “regular factors”.$^{\text{\cite{KM1}}}$%\endnote{K. Mroi, Studies in Babylonian Mathematics No.3, 2007, pp. 6-7.} 
In other words, they knew that the reciprocal of a number of the form   $ 2^{\alpha} \times  3^{\beta} \times 5^{\gamma}$    ($\alpha,\beta,\gamma$: integers) can be obtained by a finite sexagesimal expansion.

Compared with the decimal system, the existence of the regular number 3 in the sexagesimal system made it easier to calculate the result of a division whose divisor is a multiple of 3, which contributed to the development of mathematics after all. The most striking example would be $ (\frac{4}{3})^{7} \approx 7;30 $, which we can confirm in a passage of Enmetena’s foundation cone written around 2400 BCE.$^{\text{\cite{KM2}}}$%\endnote{K. Muroi, The oldest example of compound interest in Sumer: Seventh power of four-thirds, 2015, \hyperlink{https://arxiv.org/abs/1407.6246}{arXiv: 1510.00330 [math.HO]}.} 
The Sumerian scribe must have calculated as follows:
\begin{align*}
\left(\frac{4}{3}\right)^{7} &= (1;20)^{7} \\
&= (1;20)\times(1;20)\times (1;20)^{5} \\
&= (1;46,40)\times(1;20)^{5}\\
 &= \cdots\\
&= 7;29,29,32,50,22,13,20 \\
&\approx 7;30.
\end{align*}
It is impossible to imagine that he obtained the result by using decimal system in his calculations because it is difficult for him to round up $ \left(\frac{4}{3}\right)^{7}  $ to $ 7.5 $: 
\begin{align*}
\left(\frac{4}{3}\right)^{7} &\approx (1.3)^{7} = 6.2748\cdots, \\
\left(\frac{4}{3}\right)^{7} &\approx (1.33)^{7} = 7.3614\dots, \\
\left(\frac{4}{3}\right)^{7} &\approx (1.333)^{7} = 7.4784\dots. 
\end{align*}
Neither did he use a fraction such as 
\[\left(\frac{4}{3}\right)^{7} = \frac{16384}{2187} = 7+\frac{1075}{2187} = 7.491\cdots,\]
which is out of Sumerian numeral system.

The Sumerian word {\fontfamily{qpl}\selectfont a-rá-gub-ba} had been passed on to Babylonian mathematics (2000-1600 BCE) as {\fontfamily{qpl}\selectfont \textit{aragubbûm}} or {\fontfamily{qpl}\selectfont \textit{a-rá-kayyamānum}}$^{\text{\cite{SMT}}}$%\endnote{In the Susa mathematical texts. A new edition of the Susa mathematical texts in a new book (Elamite Mathematics)  is being prepared for publication by the present author and N. Heydari who is a mathematician from modern Susa.} 
with the same meaning.

\section{Number signs for area}
Since the units of area and capacity vary in different places and at different historic periods, we had better confine ourselves in this discussion to those units appearing in the   texts excavated from Fara,  the ancient city   Šuruppak   in the 26th century BCE. It is a fact, however, that the metrology of Fara  in  this historical period provided a framework for the metrology of later periods.  

The basic unit of length is {\fontfamily{qpl}\selectfont ninda-DU} ($ \approx 6m$), in which the cuneiform sign {\fontfamily{qpl}\selectfont DU} is appended probably  to distinguish it from the   {\fontfamily{qpl}\selectfont ninda}   ``bread''. In later periods, the sign {\fontfamily{qpl}\selectfont DU} was almost omitted. The basic unit  of area is squared {\fontfamily{qpl}\selectfont ninda-DU}, or {\fontfamily{qpl}\selectfont sar}, from which three larger units are derived: 
\begin{align*}
	&1~\text{{\fontfamily{qpl}\selectfont iku}} = 1,40 ~\text{{\fontfamily{qpl}\selectfont sar}}~(= 100 ~\text{{\fontfamily{qpl}\selectfont sar}}),\\
	&1~\text{{\fontfamily{qpl}\selectfont bùr}} = 18 ~\text{{\fontfamily{qpl}\selectfont iku},}\\
	& 1~\text{{\fontfamily{qpl}\selectfont šár}} = 1,0 ~\text{{\fontfamily{qpl}\selectfont bùr} (= 60 bùr)}.
\end{align*}

\begin{figure}[H]
	\centering
	\includegraphics[scale=1]{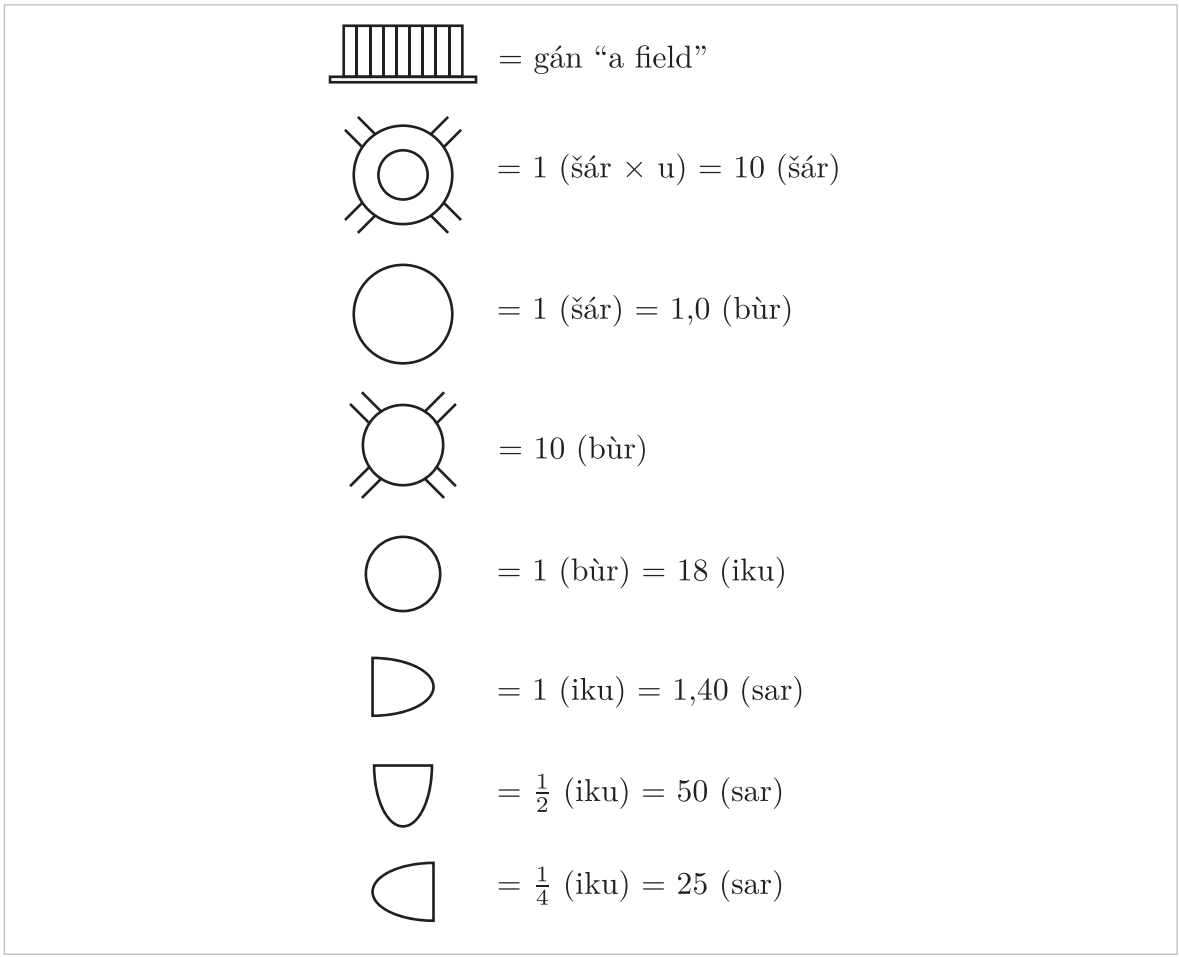}
	\caption{Number signs for area}
	\label{numbersigns}
\end{figure}

Although the number signs for area units are similar to those used in sexagesimal notation,$^{\text{\cite{MUR}}}$ we hardly mix them up if we pay close attention to the cuneiform sign  {\fontfamily{qpl}\selectfont gán} ``a field'', which is usually written down at the first line of a tablet concerning the area calculations (see Figure \ref{numbersigns}).

\section{The area of a large square field} 
\begin{figure}[H]
	\centering
	\includegraphics[scale=1]{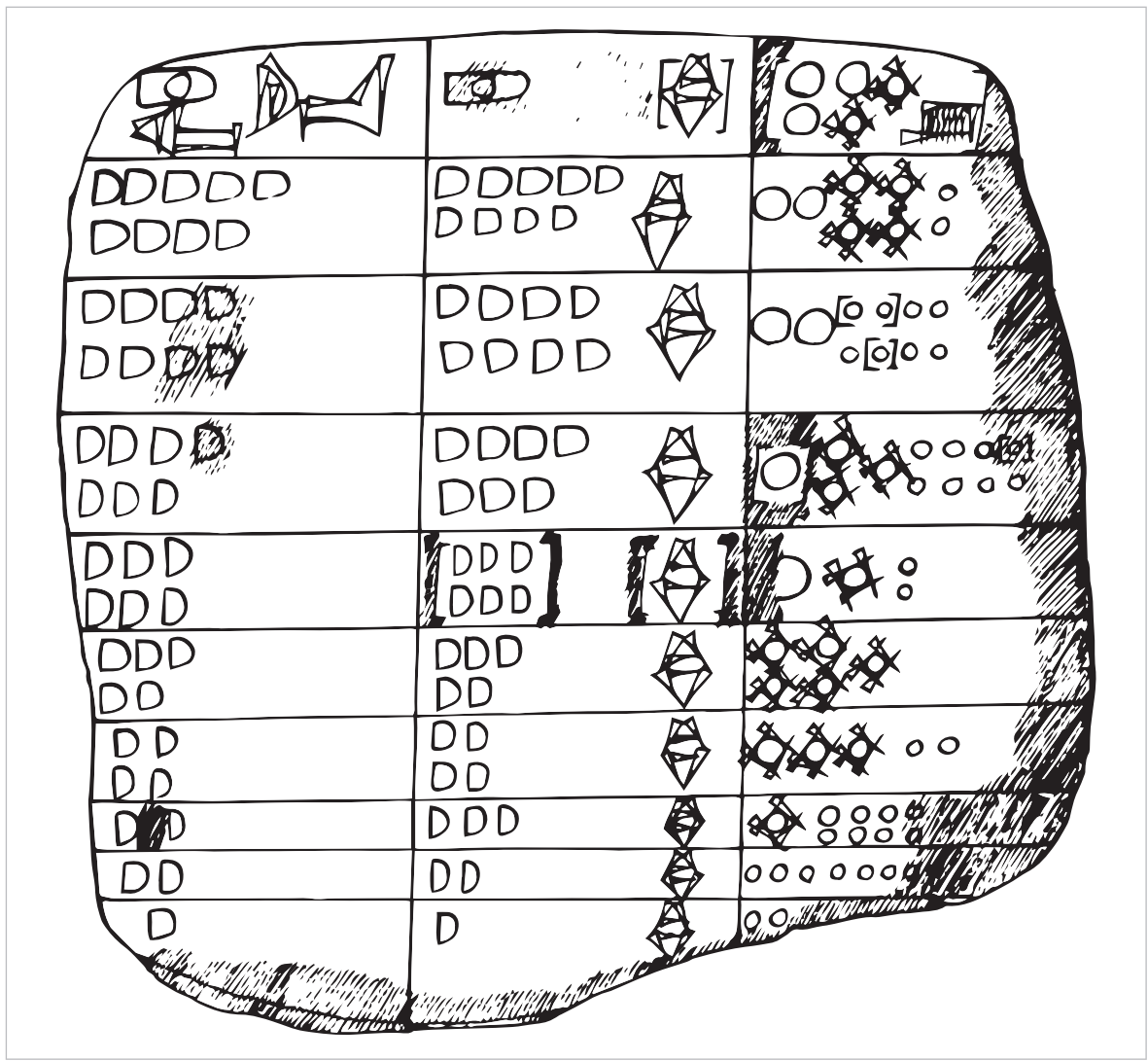}
	\caption{Obverse of School tablet \textbf{SF No.\,82}}
	\label{SF82}
\end{figure}

On the school tablet \textbf{SF No.\,82} published by A. Deimel in 1923,$^{\text{\cite{DEI1}}}$ the areas of squares whose sides decrease from 10,0 {\fontfamily{qpl}\selectfont ninda} ($ \approx 3600m  $) to 5 {\fontfamily{qpl}\selectfont ninda} ($ \approx 30m $) are tabulated, which may seem to be complicated at first sight. In order to analyze the calculation method, especially for large squares,  ten lines on the obverse of this tablet are given in the following (see Figure \ref{SF82}). Line 1, which consists of three blocks, reads as follows: 
\begin{displayquote}  
	The side is 10,0  {\fontfamily{qpl}\selectfont ninda}. (The other side) is 10,0 {\fontfamily{qpl}\selectfont ninda} equally. 3 ({\fontfamily{qpl}\selectfont šár}) 20 (bùr) is the area. 
\end{displayquote} 	

The calculation process would be as follows:
\begin{align*}
	10,0 \times 10,0~(\text{{\fontfamily{qpl}\selectfont ninda}}^2 = \text{sar})&=  1,0 \times 1,0~(\text{{\fontfamily{qpl}\selectfont éš$^2$ = iku}), where 1 {\fontfamily{qpl}\selectfont éš = 10 ninda}} \\
	&= 1,0 \times 1,0 \times \frac{1}{18} ~(\text{{\fontfamily{qpl}\selectfont bùr}), where 1 {\fontfamily{qpl}\selectfont bùr = 18 iku}} \\
	& = 3,20 ~(\text{{\fontfamily{qpl}\selectfont bùr}), where}~ \frac{1}{18} = 0;3,20 \\
	&= 3~(\text{{\fontfamily{qpl}\selectfont šár})  20 ({\fontfamily{qpl}\selectfont bùr})}.
\end{align*}

Systematically changing the area units, the Sumerian scribe has obtained the result correctly. Similarly, we can explain the following lines.\\
\underline{Line 2}:
\begin{align*}
	9,0 \times 9,0~(\text{{\fontfamily{qpl}\selectfont sar}})&=54\times 54~(\text{{\fontfamily{qpl}\selectfont iku}}) \\
	&= 54\times 54 \times \frac{1}{18} ~(\text{{\fontfamily{qpl}\selectfont bùr})} \\
	& = 2,42 ~(\text{{\fontfamily{qpl}\selectfont bùr})} \\
	&= 2~(\text{{\fontfamily{qpl}\selectfont šár}) 42 ({\fontfamily{qpl}\selectfont bùr})}.
\end{align*}  
\underline{Line 3}:
\begin{align*}
	8,0 \times 8,0~(\text{{\fontfamily{qpl}\selectfont sar}})&=48\times 48~(\text{{\fontfamily{qpl}\selectfont iku}}) \\
	&= 48\times 48 \times \frac{1}{18} ~(\text{{\fontfamily{qpl}\selectfont bùr})} \\
	& = 2,8 ~(\text{{\fontfamily{qpl}\selectfont bùr})} \\
	&= 2~(\text{{\fontfamily{qpl}\selectfont šár}) 8 ({\fontfamily{qpl}\selectfont bùr})}.
\end{align*} 
\underline{Line 4}:
\begin{align*}
	7,0 \times 7,0~(\text{{\fontfamily{qpl}\selectfont sar}})&=42\times 42~(\text{{\fontfamily{qpl}\selectfont iku}}) \\
	&= 42\times 42 \times \frac{1}{18} ~(\text{{\fontfamily{qpl}\selectfont bùr})} \\
	& = 1,38 ~(\text{{\fontfamily{qpl}\selectfont bùr})} \\
	&= 1~(\text{{\fontfamily{qpl}\selectfont šár}) 38 ({\fontfamily{qpl}\selectfont bùr})}.
\end{align*} 
\underline{Line 5}:
\begin{align*}
	6,0 \times 6,0~(\text{{\fontfamily{qpl}\selectfont sar}})&= 36\times 36~(\text{{\fontfamily{qpl}\selectfont iku}}) \\
	&= 36\times 36 \times \frac{1}{18} ~(\text{{\fontfamily{qpl}\selectfont bùr})} \\
	& = 1,12 ~(\text{{\fontfamily{qpl}\selectfont bùr})} \\
	&= 1~(\text{{\fontfamily{qpl}\selectfont šár}) 12 ({\fontfamily{qpl}\selectfont bùr})}.
\end{align*}
\underline{Line 6}:
\begin{align*}
	5,0 \times 5,0~(\text{{\fontfamily{qpl}\selectfont sar}})&= 30\times 30~(\text{{\fontfamily{qpl}\selectfont iku}}) \\
	&= 30\times 30 \times \frac{1}{18} ~(\text{{\fontfamily{qpl}\selectfont bùr})} \\
	& = 50 ~(\text{{\fontfamily{qpl}\selectfont bùr})}.
\end{align*}
\underline{Line 7}:
\begin{align*}
	4,0 \times 4,0~(\text{{\fontfamily{qpl}\selectfont sar}})&= 24\times 24~(\text{{\fontfamily{qpl}\selectfont iku}}) \\
	&= 24\times 24 \times \frac{1}{18} ~(\text{{\fontfamily{qpl}\selectfont bùr})} \\
	& = 32 ~(\text{{\fontfamily{qpl}\selectfont bùr})}.
\end{align*}
\underline{Line 8}:
\begin{align*}
	3,0 \times 3,0~(\text{{\fontfamily{qpl}\selectfont sar}})&= 18\times 18~(\text{{\fontfamily{qpl}\selectfont iku}}) \\
	&= 18\times 18 \times \frac{1}{18} ~(\text{{\fontfamily{qpl}\selectfont bùr})} \\
	& = 18 ~(\text{{\fontfamily{qpl}\selectfont bùr})}.
\end{align*}
\underline{Line 9}:
\begin{align*}
	2,0 \times 2,0~(\text{{\fontfamily{qpl}\selectfont sar}})&= 12\times 12~(\text{{\fontfamily{qpl}\selectfont iku}}) \\
	&= 12\times 12 \times \frac{1}{18} ~(\text{{\fontfamily{qpl}\selectfont bùr})} \\
	& = 8 ~(\text{{\fontfamily{qpl}\selectfont bùr})}.
\end{align*}
\underline{Line 10}:
\begin{align*}
	1,0 \times 1,0~(\text{{\fontfamily{qpl}\selectfont sar}})&= 6\times 6~(\text{{\fontfamily{qpl}\selectfont iku}}) \\
	&= 6\times 6 \times \frac{1}{18} ~(\text{{\fontfamily{qpl}\selectfont bùr})} \\
	& = 2 ~(\text{{\fontfamily{qpl}\selectfont bùr})}.
\end{align*}

For comparison, I take up another interpretation of the numbers listed in the table, which seems to have been supported by many Sumerologists:$^{\text{\cite{KRA}}}$
\begin{align*}
	(60\times 9) (60\times 9) &= 1080\times 2 + 180\times 4 + 18\times 2\\
	&= 2916 ~(\text{{\fontfamily{qpl}\selectfont iku})}.
\end{align*}
Also, see Line 2 above.

We can immediately realize that this calculation is not the one performed by the Sumerians but a modern one forced by decimal calculations, because they had no means to record the result of $(60\times 9) (60\times 9)$ decimally and their purpose was to represent the area of a large land not by unit  {\fontfamily{qpl}\selectfont iku} but by unit {\fontfamily{qpl}\selectfont bùr}. In ancient times, Sumerian framers used to sow a field with barley, a certain unit of barley ($\approx$ 240 liters) per 1 {\fontfamily{qpl}\selectfont bùr}.

\section{Incorrect result}
\begin{figure}[H]
	\centering
	\includegraphics[scale=1]{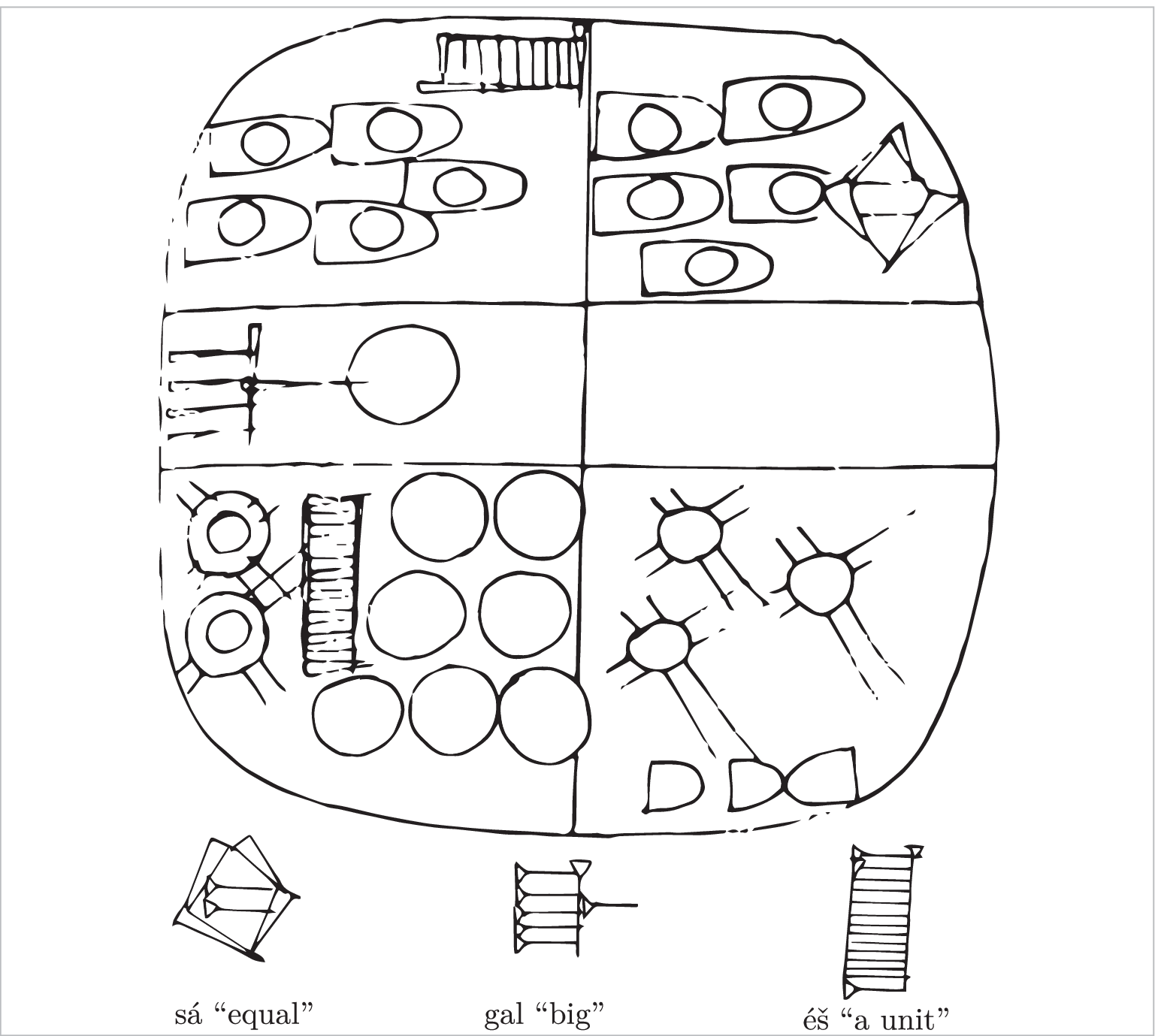}
	\caption{Obverse of School tablet \textbf{TSŠ No.\,188}}
	\label{TSS188}
\end{figure}

Another school tablet \textbf{TSŠ No.\,188}$^{\text{\cite{JES1}}}$ also calculates the area of a vast square field, whose side is definitely an imaginary one, that is, 50,0 {\fontfamily{qpl}\selectfont ninda} ($ \approx 18km$). Line 1 of the tablet consisting of two blocks simply reads (see Figure \ref{TSS188}):
\begin{displayquote}  
	~~~~~~~~~~~~~~~~~~~The field.  50,0  ({\fontfamily{qpl}\selectfont ninda}). 50,0 ({\fontfamily{qpl}\selectfont ninda}) equally.  
\end{displayquote}

If we follow the Sumerian calculation method discussed above, we can obtain the area easily:
\begin{align*}
	50,0 \times 50,0~(\text{{\fontfamily{qpl}\selectfont sar}})&= 5,0\times 5,0~(\text{{\fontfamily{qpl}\selectfont iku}}) \\
	&= 5,0\times 5,0 \times \frac{1}{18} ~(\text{{\fontfamily{qpl}\selectfont bùr})} \\
	& = 25,0,0 \times 0;3,20 ~(\text{{\fontfamily{qpl}\selectfont bùr})}\\
	& = 1,23,20 ~(\text{{\fontfamily{qpl}\selectfont bùr})}\\
	& = 1  ~(\text{{\fontfamily{qpl}\selectfont šár-gal})} 23 ~(\text{{\fontfamily{qpl}\selectfont šár})} 20 ~(\text{{\fontfamily{qpl}\selectfont bùr})}\\
\end{align*}
where 1 {\fontfamily{qpl}\selectfont šár-gal} is 1,0,0 {\fontfamily{qpl}\selectfont bùr}, and it was never used in everyday life. But the answer in lines 2 and 3 is a bit larger:
\begin{displayquote}  
	~~~~~~~~~~~~~~~1 ({\fontfamily{qpl}\selectfont šár-gal}) 27 ({\fontfamily{qpl}\selectfont šár}) 30 ({\fontfamily{qpl}\selectfont bùr}) = 1,27,30 ({\fontfamily{qpl}\selectfont bùr}).  
\end{displayquote}

In addition, the cuneiform sign  {\fontfamily{qpl}\selectfont éš} and the area number $2+\frac{1}{4}$ ({\fontfamily{qpl}\selectfont iku}) are inserted in line 3, which are seemingly unnecessary but may explain the cause of an error in calculation. The sign {\fontfamily{qpl}\selectfont éš} would emphasize that the length unit in the area unit of {\fontfamily{qpl}\selectfont iku} is not {\fontfamily{qpl}\selectfont ninda} but {\fontfamily{qpl}\selectfont éš}. As to  $2+\frac{1}{4}$ ({\fontfamily{qpl}\selectfont iku}), the scribe of this tablet must have used the relation 
\[  2+\frac{1}{4}~(\text{{\fontfamily{qpl}\selectfont iku}})= 225~(\text{{\fontfamily{qpl}\selectfont sar}}) = 15^2~(\text{{\fontfamily{qpl}\selectfont sar}})\]
to change the unit   {\fontfamily{qpl}\selectfont sar}   to  {\fontfamily{qpl}\selectfont iku}, but he has made a great mistake in changing {\fontfamily{qpl}\selectfont iku} to {\fontfamily{qpl}\selectfont bùr}:
\begin{align*}
	50,0 \times 50,0~(\text{{\fontfamily{qpl}\selectfont sar}})&= 2\frac{1}{4}\times \left(\frac{50,0}{15}\right)^2~(\text{{\fontfamily{qpl}\selectfont iku}}) \\
	&= 2\frac{1}{4}\times \left(3,20\right)^2~(\text{{\fontfamily{qpl}\selectfont iku}}) \\
	& = 2\frac{1}{4}\times \left(11,6,40\right)~(\text{{\fontfamily{qpl}\selectfont iku}})\\
	& = 2\frac{1}{4}\times \left(11,6,40\right) \times \frac{1}{18}~(\text{{\fontfamily{qpl}\selectfont bùr}})\\
	& = 2\frac{1}{4}\times \left(11,6,40\right) \times (0;3,30)~(\text{{\fontfamily{qpl}\selectfont bùr}})\\
	& = 2\frac{1}{4}\times   (38,53;20)~(\text{{\fontfamily{qpl}\selectfont bùr}})\\
	& = 1,17,46;40+9,43;20~(\text{{\fontfamily{qpl}\selectfont bùr}})\\
	& = 1,27,30~(\text{{\fontfamily{qpl}\selectfont bùr}})\\
	& = 1~(\text{{\fontfamily{qpl}\selectfont šár-gal}}) 27~(\text{{\fontfamily{qpl}\selectfont šár}}) 30~(\text{{\fontfamily{qpl}\selectfont bùr}}).
\end{align*}

In the division by 18, if the dividend is not a simple multiple of 18, the division must have been replaced by the multiplication by 0;3,20 as in later periods.

After all, the scribe's idea to use the relation  
\[  2+\frac{1}{4}~(\text{{\fontfamily{qpl}\selectfont iku}})= 225~(\text{{\fontfamily{qpl}\selectfont sar}}) = 15^2~(\text{{\fontfamily{qpl}\selectfont sar}})\]
for the conversion of units was a failure, although the idea itself was mathematically correct. 

\section{Capacity units}
\begin{figure}[H]
	\centering
	\includegraphics[scale=1]{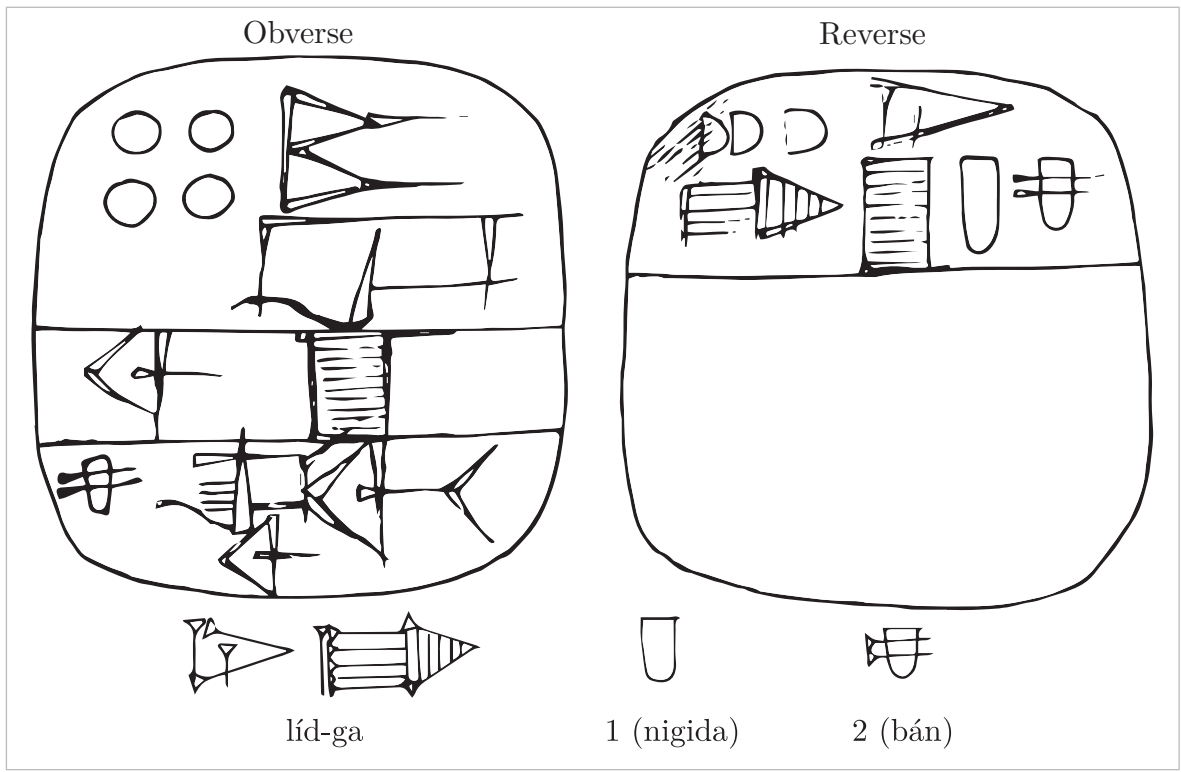}
	\caption{Obverse and reverse of administrative tablet \textbf{TSŠ No.\,81}}
	\label{TSS81}
\end{figure}

The basic capacity unit in Fara was {\fontfamily{qpl}\selectfont sìla} ($ \approx $ 1 liter), from which larger units {\fontfamily{qpl}\selectfont bán}, {\fontfamily{qpl}\selectfont nigida}, and {\fontfamily{qpl}\selectfont líd-ga} were derived:
\begin{align*}
	&1~\text{{\fontfamily{qpl}\selectfont bán}} = 10~\text{{\fontfamily{qpl}\selectfont sìla}},\\
	&1~\text{{\fontfamily{qpl}\selectfont nigida}} =  6~\text{{\fontfamily{qpl}\selectfont bán}}=  1,0~\text{{\fontfamily{qpl}\selectfont sìla}},\\
	&1~\text{{\fontfamily{qpl}\selectfont líd-ga}} =  4~\text{{\fontfamily{qpl}\selectfont nigida}}= 24~\text{{\fontfamily{qpl}\selectfont bán}}= 4,0~\text{{\fontfamily{qpl}\selectfont sìla}}.
\end{align*}
An administrative tablet of this period  \textbf{TSŠ No.\,81} confirms the  application of these units:
\begin{displayquote}  
	40 young bricklayers were hired. (Each of them) received 2 ({\fontfamily{qpl}\selectfont bán}) (of flour). (The total amount of) flour is 3 ({\fontfamily{qpl}\selectfont líd-ga}) 1 ({\fontfamily{qpl}\selectfont nigida}) 2 ({\fontfamily{qpl}\selectfont bán}).
\end{displayquote} 	
Since 
\begin{align*}
	40 \times 2~(\text{{\fontfamily{qpl}\selectfont bán}})&= 1,20~(\text{{\fontfamily{qpl}\selectfont bán}}) \\
	&= 3\times 24 + 6 + 2 ~(\text{{\fontfamily{qpl}\selectfont bán}}) \\
	& = 3  ~(\text{{\fontfamily{qpl}\selectfont líd-ga}}) 1 ~(\text{{\fontfamily{qpl}\selectfont nigida}}) 2 ~(\text{{\fontfamily{qpl}\selectfont bán}})
\end{align*}
the result given by the scribe is correct (see Figure \ref{TSS81}). Note that two words  {\fontfamily{qpl}\selectfont \textbottomtiebar{h}un} ``to hire'' and {\fontfamily{qpl}\selectfont  zì} ``flour'' are written with the same cuneiform sign  {\fontfamily{qpl}\selectfont éš} ``a unit of length''. 

There was another large unit next to the {\fontfamily{qpl}\selectfont líd-ga}, which, strangely enough, has been overlooked by scholars for almost one hundred years:
\begin{displayquote}  
	~~~~~~~~~~~~~~~~~~~~~~1  {\fontfamily{qpl}\selectfont gur}  (or the like) =  5,0   {\fontfamily{qpl}\selectfont sìla}  (= 300 {\fontfamily{qpl}\selectfont sìla}).
\end{displayquote} 	
An evidence for the use of this large unit is given in the fragmentary administrative tablet \textbf{TSŠ No.\,882}:
\begin{displayquote}  
	2  ({\fontfamily{qpl}\selectfont nigida})  of barley and 2  ({\fontfamily{qpl}\selectfont nigida}) of emmer wheat have been received as a monthly food. 1 ({\fontfamily{qpl}\selectfont nigida}) of barley and 1 ({\fontfamily{qpl}\selectfont nigida}) of emmer wheat (as a food for a half of) the ... month. $3+\frac{1}{2}$ months. In total, 1 ({\fontfamily{qpl}\selectfont gur}) 2 ({\fontfamily{qpl}\selectfont nigida}) of barley in ({\fontfamily{qpl}\selectfont gur-ma\textbottomtiebar{h}}) unit. 1  ({\fontfamily{qpl}\selectfont gur}) 2 ({\fontfamily{qpl}\selectfont nigida}) of emmer wheat. The house of gods eats (these).
\end{displayquote} 	
Thus, the total amount of barley or emmer wheat is  
\begin{align*}
	(3\times 2+1)~(\text{{\fontfamily{qpl}\selectfont nigida}})= (5+2)~(\text{{\fontfamily{qpl}\selectfont nigida}}) 
	= 1  ~(\text{{\fontfamily{qpl}\selectfont gur}})~  2 ~(\text{{\fontfamily{qpl}\selectfont nigida}}),
\end{align*}
which clearly shows that 
\begin{align*}
	1  ~(\text{{\fontfamily{qpl}\selectfont gur}})  = 5 ~(\text{{\fontfamily{qpl}\selectfont nigida}}) = 5,0~(\text{{\fontfamily{qpl}\selectfont sìla}}).
\end{align*}
For more attestations of this capacity unit, see \textbf{WF No.\,1}  and  \textbf{WF No.\,7}$^{\text{\cite{DEI2}}}$ in the following section.

\section{Barley distributed to donkeys}
There exists a group of tablets excavated in Fara in which the total amount of barely distributed to {\fontfamily{qpl}\selectfont anše-apin} ``plow-donkeys'' is recorded with the number of donkeys and the names of their owners. It seems that 1 ({\fontfamily{qpl}\selectfont nigida}) of barely was assigned to one donkey probably for sowing a field, and therefore the formula for calculating the total amount of barley would be:
\begin{displayquote}  
	the number of donkeys divided by 4 in the case of 	{\fontfamily{qpl}\selectfont líd-ga} (= 4 {\fontfamily{qpl}\selectfont nigida}), or the number of donkeys divided by 5 in the case of {\fontfamily{qpl}\selectfont gur} (= 5 {\fontfamily{qpl}\selectfont nigida}).
\end{displayquote} 	

In the following, several examples  of these types from \textbf{WF}, \textbf{TSŠ} and \textbf{NTSŠ}$^{\text{\cite{JES2}}}$  are given:\\

\noindent
\textbf{WF No.\,1}:
\[ \mathbf{1,26}\div 5= 17 +\frac{1}{5}~(\text{{\fontfamily{qpl}\selectfont gur}}) = 17 ~(\text{{\fontfamily{qpl}\selectfont gur}})~ 1 ~(\text{{\fontfamily{qpl}\selectfont nigida}})  \]
where the number in  boldface   represents the head of plow-donkeys, as above. \\

\noindent
\textbf{WF No.\,7}:
\[ \mathbf{6,23}\div 5= 1,16 +\frac{3}{5}~(\text{{\fontfamily{qpl}\selectfont gur}}) = 1,16 ~(\text{{\fontfamily{qpl}\selectfont gur}})~ 3 ~(\text{{\fontfamily{qpl}\selectfont nigida}}).  \]

\noindent
\textbf{WF No.\,14}:
\[ \mathbf{11}\div 4= 2 +\frac{3}{4}~(\text{{\fontfamily{qpl}\selectfont líd-ga}}) = 2 ~(\text{{\fontfamily{qpl}\selectfont líd-ga}})~ 3 ~(\text{{\fontfamily{qpl}\selectfont nigida}}).  \]

\noindent
\textbf{TSŠ No.\,115}:
\[ \mathbf{1,17}\div 4= 20 - \frac{3}{4}~(\text{{\fontfamily{qpl}\selectfont líd-ga}}) = 20 ~(\text{{\fontfamily{qpl}\selectfont líd-ga}})- 3 ~(\text{{\fontfamily{qpl}\selectfont nigida}}).  \]

\noindent
\textbf{NTSŠ No.\,211}:
\[ \mathbf{42}\div 4= 10 +  \frac{2}{4}~(\text{{\fontfamily{qpl}\selectfont líd-ga}}) = 10 ~(\text{{\fontfamily{qpl}\selectfont líd-ga}}) ~2 ~(\text{{\fontfamily{qpl}\selectfont nigida}}).  \]

These formulas are very simple and practical, however it is really regrettable that most Sumerologists have not  noticed their existence yet.

\section{Conclusion}
Up to the present, most Sumerologists have transcribed the Sumerian sexagesimal numbers to our decimal numbers, both in papers and in books,$^{\text{\cite{SNK}}}$%\endnote{For example, S. N. Kramer, The Sumerians, 1963, p. 92.}
 misleading the general readers and perhaps themselves to the assumption that the Sumerians calculated decimally and wrote down only the results in the sexagesimal numeral notation. It is, however, against the facts as we have seen above. From the beginning they calculated in the sexagesimal numeral system and recorded the results according to their weights and measures which were rather complicated.

Moreover, the Sumerians had established the reasonable and convenient systems of area and capacity measures based on their sexagesimal numeral system, which seem to have made it possible to manage farm work efficiently. By analyzing  the data recorded by the Sumerians, we can clearly recognize how they had calculated a large area of a field and a large amount of barley necessary for sowing a field. However, we must be very cautious about  transcribing the sexagesimal numbers occurred in their documents to our decimal numbers, because such a transcription tends to obscure the calculation process involved and may mislead us about the nature of Sumerians' mathematics. 

I would like to admire the calculation ability of the Sumerian scribes that must have been acquired after long and strict education.


\begin{thebibliography}{10}
\bibitem{SMT}
In the Susa mathematical texts. A new edition of the Susa mathematical texts in a new book (Elamite Mathematics)  is being prepared for publication by the present author and N. Heydari who is a mathematician from modern Susa.
\end{thebibliography}
\end{document}